\documentclass[11pt]{article}
\usepackage{amsmath}
\usepackage{amssymb}
\usepackage{amsfonts}
\usepackage{amscd}
\usepackage{epsfig}
\newtheorem{theorem}{\bf Theorem}[section]

\newtheorem{corollary}{\bf Corollary}[section]
\newtheorem{remark}{\bf Remark}[section]
\newtheorem{lemma}{\bf Lemma}[section]
\begin{document}

\def\ds{\displaystyle}
\title{
Negative norm estimates and superconvergence  results  in  Galerkin method  for  
strongly nonlinear parabolic problems} 
\author{ Ambit Kumar Pany \footnote {Center of Applied Mathematics and Computation, SOA (Deemed to be University),
Bhubaneswar-752030, India. Email: ambit.pany@gmail.com, ambitpany@soa.ac.in}\;,
Morrakot Khebchareon
\footnote {Department of Mathematics, Faculty of Science, Chiang Mai
University, Chiang Mai-50200 (Thailand) and Advanced Research Centre for Computational Simulation, Chiang Mai  University, Chiang Mai-50200 (Thailand) and   Centre of Excellence in Mathematics, CHE, 328 Si Ayutthaya Road, Bangkok,
(Thailand). Email:morrakot.k@cmu.ac.th}\;,
and  Amiya K. Pani \footnote {Department of Mathematics,
Indian Institute of Technology Bombay, Powai, Mumbai-400076 (India).
Email: akp@math.iitb.ac.in}}

\maketitle
\begin{abstract}
The conforming finite element Galerkin  method is applied to discretise in the spatial direction for a class of strongly nonlinear   parabolic problems.
Using elliptic projection of the associated  linearised stationary  problem with Gronwall type result, optimal error estimates are derived, 
when piecewise polynomials of degree $r\geq 1$ are used,
which improve upon earlier results of  Axelsson  [Numer. Math. 28 (1977), pp. 1-14] requiring for 2d $r\geq 2$ and for 3d $r\geq 3.$
Based on quasi-projection technique introduced  by Douglas {\it et al.}  [Math. Comp.32 (1978),pp. 345-362], superconvergence result for the error between Galerkin approximation and approximation through quasi-projection  is established for  the semidiscrete Galerkin scheme.  Further, {\it a priori} error estimates in Sobolev spaces of  negative index are derived. Moreover, in a single space variable, nodal superconvergence results between the true solution and Galerkin approximation are established. 


\end{abstract}
\section{Introduction}
In this paper, {\it a priori} error estimates with  superconvergence results and negative norm  estimates are derived for the conforming finite element Galerkin method applied to the following class of strongly nonlinear  parabolic initial and boundary value problems:
Find $u=u(x,t),\;x\in \Omega,\;t\in J=(0,T]$ for any $T >0$ satisfying
\begin{eqnarray}\label{eq:1.1}
\frac{\partial u}{\partial t}-\nabla\cdot A(u,{\nabla}u) + f(u,{\nabla}u) = 0, \;\;(x,t) \in \; \Omega \times J
\end{eqnarray}
with initial and homogeneous Dirichlet boundary conditions
\begin{eqnarray}
u(x,0)=u_0(x),\;\;  x \in \Omega, \label{eq:1.2}\\
u(x,t) = 0, \;\; x\in {\partial{\Omega}},\;\;\; t \in \;J,\label{eq:1.3}
\end{eqnarray}
where ${\Omega}$ is a bounded domain in  $ \mathbb{R}^d,\;d=1,2,3$  with smooth boundary $\partial \Omega.$  Further,
$ A=(A_1,\cdots,A_d) $ and $ f $ are   vector and scalar valued functions defined on $ \mathbb{R} \times
\mathbb{R}^d,$ respectively.

Throughout this article, we  make the following assumptions called  $({\bf A}_1)$ on the coefficients:
\begin{itemize}
\item [(i)] The problem (\ref{eq:1.1})-(\ref{eq:1.3}) has a unique solution with bounded gradients $ |\nabla u|, |\nabla u_t|$.
\item [(ii)] The function $A$ and $f$ are sufficiently smooth and bounded.
\item [(iii)] The matrix $\frac{\partial A}{\partial \xi} = [\frac{\partial A_i}{\partial \xi_j}]$, where $\xi =(\xi_1,\xi_2,\ldots,\xi_d)$ and $\xi_j = \frac{\partial u}{\partial x_j}$ is uniformly positive definite.
\end{itemize}
For existence, uniqueness and regularity results of such nonlinear equations, see \cite{LSU-1968}. For related numerical methods,  a good number of article is devoted to strongly nonlinear elliptic problems, see, \cite{BL-2012}, \cite{C-1989}, \cite{FZ-1987}, \cite{GNP-2008}, \cite{HRS-2005}, \cite{MP-1995}, \cite{P-1995}, \cite{P-2005}, \cite{SP-2018} and references, there in. However, there seems to be less number of papers available in literature on numerical approximation to strongly nonlinear parabolic problems, see, \cite{A-1977}, \cite{KMP-1996}, \cite{SKP-2020} and \cite{KPK-2016}, etc. The more relevant article is  \cite{A-1977}, where conforming FEM is applied to the problem \eqref{eq:1.1}-\eqref{eq:1.3} and optimal error estimates in $L^{\infty}(L^2)$  are derived using piecewise polynomial of degree $r\geq 2$ for $d=2$ and for $d=3$, $r\geq 3.$ One of our effort is to prove  global optimal error estimates when $r\geq 1.$

Superconvergece, one of the objectives  of  this article, has long been an active area of  research in finite element methods, see \cite{LW-1995}, \cite{KN-1987}, \cite{KNS-1997} and \cite{CZ-2018}.This is  mainly due to its applications in {\it a posteriori } error estimates. One prominent line of thought is to post-process the discrete solution. The major tool used in this paper is an asymptotic expansion called quasi-projection of the approximate solution which leads to an expansion of the error.  Essentially, being a postprocessing method,  it is based on a sequence of projections associated with the approximate solution of the underlying stationary problem.
Earlier, Douglas {\it et al.} \cite{DDW-1978} introduced  quasi-projection technique  for linear parabolic and second order hyperbolic equations   and analyzed the superconvergence phenomena associated with  1D-problems with the help of  negative norm estimates. Later, Arnold and Douglas \cite {AD-1979} generalized  these results  to a quasilinear parabolic equation of the form
$ c(x,t,u)u_t - \nabla \cdot [a(x,t,u)\nabla u + b(x,t,u)]+f(x,u,t)=0.$ For similar results for parabolic integro-differential equations, see, \cite{PS-1996} and for the Stefan problem, see, \cite{LP-1995}.
For detailed survey of  superconvergence result,  refer to   \cite{KN-1987}, \cite{KNS-1997} and  \cite{CZ-2018} and  references, there in.
Regarding the application of negative norm estimates to  prove interior superconvergence results is not new to the literature, see, \cite{BS-1977}\cite{BSTW-1977}, and \cite{T-1980}. We strongly believe that the tools for global negative norm estimates will help to prove interior estimates and using the postprocessing technique  in terms of certain averaging of operator applied to the Galerkin approximations,see \cite{BS-1977}, it is possible to prove interior superconvergence results, which will be a part of our future  work. 


The main contributions  of this paper are summerised  as follows.
\begin{itemize}
\item Optimal  error estimates of the semidiscrete Galerkin approximation are derived for the problem \eqref{eq:1.1} in  $d=1,2,3.$ 

Compared to  Axelsson \cite{A-1977},  additional superconvergence estimate in $L^{\infty}(H^1)$  for the error between the Galerkin approximation and  the elliptic projection  is proved and as a consequence, optimal estimates in $L^{\infty}(H^1)$ norm for $d=1,2,3$ and   in $L^{\infty}(L^{\infty})$ -norm when $d=1,2$ are shown.  The optimal error estimate in $L^{\infty}(L^2)$ is valid for $d=1,2,3$, when piecewise polynomials of degree $r\geq 1$ are used. But in \cite{A-1977},  optimal estimate is valid for  $r\geq 2,$ that is, for quadratic or higher order elements in a one or two dimensional problems, where as in a three  dimensional problems,   optimal bound is derived
 for  $r\geq 3,$ that is, for cubic or higher order elements.

\item Based on quasi-projection technique of Douglas {\it et al.}   \cite{DDW-1978},  superconvergence results are  established  for the error between Galerkin approximation and quasi-projection. 
As a result, optimal negative norm estimates are shown for the error between the semidiscrete Galerkin approximation and the exact solution.

\item In  a single space variable, knot superconvergence of  semidiscrete Galerkin approximation to the solution of (\ref{eq:1.1})-(\ref{eq:1.3}) is derived  using quasi-projection technique of \cite{DDW-1978}.
\end{itemize}

A general outline of the paper is as follows: In section 3, some basic result are derived. A quasi-projection is defined section 4 and section 5 deals with superconvergence phenomenon for a single space variable. The convergence at the knot points are shown to be of order $h^{2r-\frac{1}{2}}$. Some of the results which are assumed in section 5 are proved in section 6. In addition, an optimal $L^\infty$ estimate is also derived in this section.

\section{Some Notations and Preliminaries.}
\setcounter{equation}{0}
In this section, we first recall the usual definitions of
standard  Sobolev spaces $W^{m,p}(\Omega)$ with the norm
$$
\|u\|_{m,p}:= (\sum_{|\alpha|\leq m}\|D^{\alpha}u\|^p \, dx)^{1/p},
\, 1\leq p <\infty,
$$
and for $p=\infty$
$$
\|u\|_{m,\infty}:= \max_{|\alpha|\leq m} \|D^{\alpha}u\|_{L^\infty (\Omega)}.
$$

When $p=2,$ denote $W^{m,2}(\Omega)$ by $H^m(\Omega)$ with norm $\|\cdot\|_m = \|\cdot\|_{H^m(\Omega)}.$
If $m=0,$ then set $\|\cdot\|=\|\cdot\|_{L^2 (\Omega)}.$ Further,
let $H^1_0 = \{ v\in H^1(\Omega): v=0 \;\;\;{\mbox{ in }}\; \partial\Omega \}$ and
let $H^{-s}$ be the dual space of $H^s$ with the norm
$$
\|v\|_{-s}=\sup_{{ w\in H^s}\,{\|w\|_s \neq 0}} \frac {|<v, w>|}{\|w\|_s}.
$$
For a normed linear space $X$ with norm $\|\cdot\|_X $, let
$$
L^p(0,T;X)=\{\phi : (0,T]\rightarrow X \,\,\, with\,\,\|\phi\|_{L^p(0,T;X)}=
\left[\ds\int^T_0\|\phi(t)\|^p_X \;dt\right]^{\frac{1}{p}}< \infty\}
$$ and
$$
L^{\infty}(0,T; X)=\{\phi : [0,T]\rightarrow  X\;\;| \;\;\|\phi\|_{L^{\infty}(0,T; X)}=
{ \rm ess sup}_{0\leq x \leq T} \|\phi(t)\|_X \}.
$$

The weak formulation of the problem (\ref{eq:1.1})-(\ref{eq:1.3}) is to find $u(t)\in H^1_0(\Omega)$ for
$t\in (0,T]$ such that
\begin{eqnarray} \label{eq:weak-formulation}
(u_t,v)+(A(u,\nabla u),\nabla v)+(f(u,\nabla u),v)= 0,\;\;v \in H^1_0(\Omega).
\end{eqnarray}

Throughout this paper, $C$ denotes a generic positive constant.
By $C(q,k)$, we mean a generic constant depending on$\|\frac {\partial^j u}{\partial t^j}\|_{L^\infty(W^{q,\infty})}$, for $j = 0,1,2,\ldots,k$, but on no higher derivatives of $u.$  For simplicity, we write $u_t$ for $\frac {\partial u}{\partial t},$ $A_u$ for $\frac{\partial A}{\partial u}$ and $A_\xi$ for $\frac {\partial A}{\partial \xi}$.  Similarly for $f_{\xi}$ and $f_{u}$.

\section{Semidiscrete Galerkin Approximation}
\setcounter{equation}{0}
This section deals with the semidiscrete formulation and elliptic projection with related error estimates.

Let  $ \mathcal{T}_h $ be  a family of shape regular and quasi-uniform  triangulation of
$\bar{\Omega}$ into simplexes $K.$ Let the discretization
parameter $h$ be defined as $h = \ds \max_{K \in \mathcal{T}_h}
{h_{K}}$, where $h_K$ is the diameter of $ K.$
Further, let $\{V^0_h\}$ for $0 < h \leq 1$ be a family of finite element  subspace of $H_0^1$
 defined by
$$V^0_h=\{v \in
C^0({\bar \Omega})\cap H^1_0:~~ v|_K \in P_r(K) \quad \forall K \in
\mathcal{T}_h\},$$
where $P_r(K)$ is space of polynomials of degree less than equal to $r$.
Note that  $V^0_h$ for $0 < h \leq 1$   satisfies the following approximation property:
There exists a constant $C$ independent of $h$ such that for $\phi \in H^q(\Omega)\cap H^1_0 $ with  $ q \in [1,r+1]$
\begin{eqnarray}\label{eq:approx-property}
\inf_{\chi \in V^0_h} \{\|\phi - \chi\| + h\|\phi - \chi\|_1\} \leq  C h^q \|\phi\|_q.
\end{eqnarray}
In addition to (2.1), the following inverse assumption on $V_h^0$  holds
for $v_h \in V_h^0$ and  for  $K\in \mathcal{T}_h $
\begin{eqnarray}\label{eq:inverse}
\|v_h\|_{H^1(K)} \leq C h_{K}^{-1}\|v_h\|_{L^2(K)} \;\;\;\mbox{and}\;\; \|v_h\|_{W^{j,\infty} (K)}
\leq C h_{K}^{-\frac{d}{2}}\|v\|_{H^j(K)}\;\;j=0,1.
\end{eqnarray}
Note that for $1\leq p\leq q\leq \infty$ and $\chi \in {V}_h^0$, the following  property  is valid:
\begin{equation}\label{eq:inverse-1}
\|\chi\|_{L^p(K)} \leq C  h_{K}^{d(\frac{1}{p}-\frac{1}{q})}\; \|\chi\|_{L^q(K)} \quad \forall K\in  \mathcal{T}_h.
\end{equation}
For properties \eqref{eq:approx-property}-\eqref{eq:inverse-1}, see \cite{BS-2008}.

The semidiscrete Galerkin approximation is now defined as a solution $U(t) \in V_h^0$  for $t\in J$ of
\begin{eqnarray}\label{eq:semi-discrete}
(U_t,V)+(A(U,\nabla U),\nabla V)+(f(U,\nabla U),V)= 0,\;V \in V^0_h
\end{eqnarray}
with $U(0)=U_0\in V^0_h$  to be defined appropriately later on as an approximation of $u_0.$

For $u,v,W\in H^1_0(\Omega),$ define
\begin{eqnarray} \label{eq:3.3}
a(u;W,v)= (A(W,\nabla W)-A(u,\nabla u),\nabla v) +(f(W,\nabla W)-f(u,\nabla u),v).
\end{eqnarray}
As in \cite{A-1977}, it is easy to check that
\begin{eqnarray}
a(u,W;W-u) \geq \rho\|\nabla(W-u)\|^2-\rho_0\|(W-u)\|^2,\;\;\;u,W \in H^1_0(\Omega),
\end{eqnarray}
where $ \rho = \inf_{\{W,\nabla W \in {\mathbb R}\times{\mathbb R}^n\}}$ \{smallest eigenvalues of $\frac{\partial A}{\partial \xi}(W,\nabla W)$\} $> 0$ and the constant $\rho_0 = \sup_{W,\nabla W}\Big\{\frac{1}{2} div [\frac{\partial A(W,\nabla W)}{\partial u} + \frac{\partial f(W,\nabla W)}{\partial u} ]\Big\}$.

Now using Taylor's expansion as in \cite{A-1977}, we write
\begin{eqnarray}\nonumber
&&a(u,W;V)=b(u,\nabla u;W-u,V) \\
&&\;\;\;\;\;+\int_\Omega[\nabla V^T,V]\int_0^1(1-s){\mathbb A}^{'} (w,\nabla w)ds
\left [\begin{array}{c}
\nabla(W-u)\\
W-u\end{array}\right ]
\left [\begin{array}{c}
\nabla(W-u)\\
W-u\end{array}\right ]\;dx,\nonumber\\
&& \;\;\;= b(u,\nabla u;W-u,V) + ( {\mathbb R}_1(W-u,\nabla(W-u)), \nabla V) +   ( {\mathbb R}_2(W-u,\nabla(W-u)),  V), \label{3.5}
\end{eqnarray}
  where $w = u+s(W-u),\;b(u,\nabla u;\phi,\psi)= \int_\Omega [\nabla \psi^T,\psi]{\mathbb A } (u,\nabla u) \left
[\begin{array}{c}
\nabla \phi\\
\phi\end{array}\right]\;dx$,
$$
{\mathbb A} (w,\nabla w)=\left [\begin{array}{cc}
\frac{\partial A(w,\nabla w)}{\partial \xi}, & \frac{\partial A(w,\nabla w)}{\partial u}\\
\frac{\partial f(w,\nabla w)}{\partial \xi}, &\frac{\partial f(w,\nabla w)}{\partial u}\end{array}\right ],
$$
and ${\mathbb A}^{'}$ is the Fr\'{e}chet derivative of ${\mathbb A}$. For similar expressions on reminder, one refers to \cite{MP-1995} and \cite{P-1995}. The bilinear form $ B(u,\nabla u;\phi,\psi)$  associated with elliptic operator ${\mathbb L}$ is given by
$$
B(u,\nabla u;\phi,\psi):= b(u,\nabla u;\phi,\psi) + \lambda (\phi,\psi)= ({\mathbb L}(u)\phi,\psi),
$$
where
\begin{eqnarray}\nonumber
{\mathbb L}(u)\phi &=& -div(\frac{\partial A}{\partial \xi}(u,\nabla u)\nabla \phi) +[-\frac{\partial A}{\partial u}(u,\nabla u)+\frac{\partial f}{\partial \xi}(u,\nabla u)]^T\nabla \phi \\
&&\;\;\;\;+ [-div\frac{\partial A}{\partial u}(u,\nabla u)+\frac{\partial f}{\partial u}(u,\nabla u)]\phi + \lambda \phi.
\end{eqnarray}
Here, $\lambda>0$ to be chosen large so that the bilinear form $B(u,\nabla u; \cdot,\cdot)$ is coercive in the sense that
$$B(\cdot,\cdot; v,v)\ge \rho\|\nabla v\|^2 + (\lambda - C(2,0))\;\|v\|^2 \geq \alpha_0 \|\nabla v\|^2.$$
Note  that $B(u,\nabla u; \cdot,\cdot)$ is also bounded, that is, $|B(u,\nabla u;v,w)| \le C \|\nabla v\|\;\|\nabla w\|.$
With $a_{\lambda}(u,W;v)= a(u,W;v) +\lambda (W-u, v),$  it satisfies the following coercivity:
$$a_{\lambda} (u,W; W-u)\ge \alpha_0 \|\nabla (W-u)\|^2.$$

The Dirichlet  problem for the bilinear form $B(u,\nabla u;\cdot,\cdot)$ has a unique solution, see \cite{DDS-1971}.
Let ${\mathbb L}^*$ be the adjoint of ${\mathbb L}$ and $\psi \in H^{s}$ satisfies
\begin{eqnarray}\label{eq:adjoint-1}
{\mathbb L}^*(u)\phi&=& \psi,\;\;x \in \Omega,\\
\phi &=&  0 , \;\;x \in \partial\Omega, \label{eq:adjoint-2}
\end{eqnarray}
satisfying the following elliptic regularity, see \cite{A-1965}
\begin{eqnarray}\label{eq:elliptic-regularity}
\|\phi\|_{s+2} \le C(s+2,0)\|\psi\|_s.
\end{eqnarray}
    Let $\tilde {u}_h : [0,T]\rightarrow V^0_h$ be the elliptic projection of $ u $ defined by
\begin{eqnarray}\label{eq:elliptic-projection}
B(u,\nabla u;u-\tilde {u}_h,v)= 0,\;\; v \in V^0_h.
\end{eqnarray}
For a given $u$, an application of the Lax-Milgram Lemma implies
the existence of a unique $\tilde{u}_{h}.$

Let $ \zeta = \tilde u_h - U,$ $ \eta = \tilde u_h-u$ and $e = U-u = \eta - \zeta$.
Subtracting (\ref{eq:weak-formulation}) from (\ref{eq:semi-discrete}), and applying (\ref{3.5}) with  (\ref{eq:elliptic-projection}), we now arrive at 
\begin{eqnarray}
(\zeta_t,v)&+&B(u,\nabla u;\zeta,v)
=(\eta_t,v)  + \lambda (\zeta,v) -\lambda (\eta,v) \nonumber\\
&&\;\;- \int_\Omega[\nabla v^T,v]\int_0^1(1-s){\mathbb A}^{'}(w,\nabla w)ds
\left [\begin{array}{c}
\nabla e\\
e \end{array}\right ]
\left [\begin{array}{c}
\nabla e\\
e\end{array}\right ]\;dx\nonumber\\
&& \;\;=(\eta_t,v)  + \lambda (\zeta,v) -\lambda (\eta,v) 
-({\mathbb R}_1(e,\nabla e),\nabla v) +  ({\mathbb R}_2(e,\nabla e), v),
 \label{eq:zeta}
\end{eqnarray}
where $w = u + s(U-u)= u+s\;e$.

\subsection{Elliptic projection.}
This subsection deals with an elliptic projection as defined by Wheeler \cite{W-1973} and related error analysis.

We  need for our subsequent use, the  interaction of  $B(u,\nabla u;\cdot,\cdot)$  with time differentiation, when the coefficient of ${\mathbb L}$  associated with this bilinear form are time dependent.
Now, let $\phi : J\rightarrow H_0^1$ and $\psi \in {H^1_0\cap H^{s+2}}$. Then, a simple use of the  Leibnitz's rule yields
\begin{eqnarray}\label{eq:B-time-derivative}
\frac{d^k}{dt^k} B(u,\nabla u;\phi, \psi)&=& \sum^k_{i=0}\left (\begin{array}{c}
k\\ i \end{array}\right )\int_\Omega  [\nabla \psi^T,\psi][\frac{d^{k-i}}{dt^{k-i}}{\mathbb A}(u,\nabla u)]\frac{d^i}{dt^i}\left [\begin{array}{c}
\nabla \phi\\ \phi \end{array}\right ]\;dx\nonumber\\
&&\;\;\;\;\;+ \lambda ( \frac{\partial^k \phi}{\partial t^k},\psi)\nonumber\\
&=&B(u,\nabla u;\frac{\partial^k \phi}{\partial t^k},\psi)+\sum^{k-1}_{i=0}F_{ik}(u,\nabla u; \frac{\partial^i \phi}{\partial t^i},\psi),
\end{eqnarray}
where for $\gamma \in H^1_0,$
\begin{eqnarray*}
F_{ik}(u,\nabla u;\gamma,\psi)&=&\left (\begin{array}{c}k\\ i \end{array}\right ){\int_\Omega  [\nabla \psi^T,\psi][\frac{d^{k-i}}{dt^{k-i}}{\mathbb A}(u,\nabla u)] \left |\begin{array}{c} \nabla \gamma\\ \gamma \end{array}\right |dx} .\\
&=& \left (\begin{array}{c}k\\ i \end{array}\right )\left(\gamma,-\nabla (\frac{d^{k-i}}{dt^{k-i}})(\frac{\partial A}{\partial \xi})\nabla \psi
 + (\frac{d^{k-i}}{dt^{k-i}})(\frac{\partial f}{\partial \xi})\psi \right) \\
&+&(\frac{d^{k-i}}{dt^{k-i}})(\frac{\partial A}{\partial u})\nabla \psi^T + (\frac{d^{k-i}}{dt^{k-i}})(\frac{\partial f}{\partial u})\psi.
\end{eqnarray*}
Therefore,
\begin{eqnarray}\label{eq:estimate-F}
|F_{ik}(u,\nabla u;\phi,\psi)| = \left\{\begin{array} {lll}
 \displaystyle C(1,k-i)\|\phi\|_1\|\psi\|_1,  \\
\displaystyle C(s+2,k-i)\|\phi\|_{-s}\|\psi\|_{s+2},\;\;\; s=0,1,2,\ldots.
\end{array} \right.
\end{eqnarray}
For convenience, we prove the following lemma, we refer \cite{AD-1979}  for an analogous result.
\begin{lemma}\label{lemma:3.1}
Let there be given a linear functional $ F :H^1_0 \rightarrow {\mathbb R}$ and numbers $ M_1\ge M_2\ge M_3 \ge \ldots \ge M_{p+1},0 \le p\le r $ with
$$
|F(\rho)| \le M_{s+2}\|\rho\|_{s+2}\;\; \mbox{for} \;\;\rho \in H^{s+2}(\Omega)\cap H^1_0(\Omega), s= -1,0,1,\ldots,p-1.
$$
Suppose $\varphi \in H^1_0(\Omega)$ satisfies
\begin{eqnarray}
B(u,\nabla u;\varphi,\chi)= F(\chi) \;\;\; \mbox{for} \quad \chi \in V^0_h,
\end{eqnarray}
then, 
\begin{eqnarray}\nonumber
\|\varphi\|_{-s} \le C(max(s,0)+2,0)\left[(M_1+\inf_{\chi \in V^0_h}\|\varphi - \chi\|_1)h^{s+1}+M_{s+2}\right] \\
 s= -1,0,1,\ldots,p-1.
\end{eqnarray}

\end{lemma}
\noindent
{\it Proof}. For $s=-1$, we note that
\begin{eqnarray*}
B(u,\nabla u;\varphi,\varphi)&=&[B(u,\nabla u;\varphi,\varphi-\chi)+F(\chi-\varphi)+F(\varphi)],\quad \chi \in M^0_h \\
&\le& \Big( C \|\varphi\|_1+M_1\Big)\;\inf_{\chi \in V^0_h} \|\varphi - \chi\|_1 + M_1\|\varphi\|_1.
\end{eqnarray*}
Since $B(u,\nabla u;\cdot,\cdot)$ is coercive, it now follows that
$$
\|\varphi\|_1\le C(2,0)\left(M_1 +\inf_{\chi \in V^0_h}\|\varphi - \chi\|\right),
$$
which yields the desired for $s = -1$.

For $s=0,1,2,\ldots,p-1,$  we apply Aubin-Nitsche's duality argument, Given $\psi \in H^s(\Omega),\;0\le s \le p-1,$ define $\phi \in H^{s+2}(\Omega)\cap H^1_0(\Omega)$ by (\ref{eq:adjoint-1})-(\ref{eq:adjoint-2}).
Now, a use of boundedness of the bilinear form with  bound for $F$, approximation property (\ref{eq:approx-property})  and (\ref{eq:elliptic-regularity}) yields
\begin{eqnarray}\nonumber
(\varphi,\psi) &=& (\varphi,{\mathbb L}^{*}(u)\phi)=B(u,\nabla u;\varphi,\phi)\\ \nonumber
&=&B(u,\nabla u;\varphi,\phi-\chi)+F(\chi - \phi)+F(\phi),\;\;\;\chi \in V^0_h \\ \nonumber
&\le& C(s+2,0)(\|\varphi\|_1 + M_1)\inf_{\chi \in V^0_h}\|\phi - \chi\|_1+M_{s+2}\|\phi\|_{s+2}\\ \nonumber
&\le& C(s+2,0)(\|\varphi\|_1 + M_1)h^{s+1}\|\phi\|_{s+2}+M_{s+2}\|\phi\|_{s+2} \\
&\le& C(s+2,0)([\|\varphi\|_1 + M_1]h^{s+1}+M_{s+2})\|\psi\|_{s}.\;\;\; 
\end{eqnarray}
For $s=0 $, we  obtain
$$
\|\varphi\|\le C(2,0)[(\|\varphi\|_1 + M_1)h+M_2)],
$$
which completes  the desired  result for $s = 0$. The other negative estimate follows similarly and this concludes the rest of the proof. 
\hfill{$\Box$}

Below, we discuss the negative norm estimate for $\eta_t$.
\begin{theorem}\label{theorem:eta-t}
Let $ 1 \le q \le {r+1}$ and $ \frac{\partial^k u}{\partial t^k} \in H^q(\Omega)$, for $ t \in J.$
Then, it holds
\begin{eqnarray}\label{eq:eta-negative}
\|\frac{\partial^k \eta}{\partial t^k}\|_{-s} \le C(max(q,s+2),k) h^{s+q}, -1 \le s \le r-1.
\end{eqnarray}
\end{theorem}
\noindent
{\it  Proof}.
A use of (\ref{eq:B-time-derivative}) shows
$$
B(u,\nabla u;\frac{\partial^k \eta}{\partial t^k},\chi)=-\sum^{k-1}_{i=0} F_{ik}(u,\nabla u;\frac{\partial^i \eta}{\partial t^i},\chi),\;\;\; \chi \in V^0_h.
$$
Now, we identify the right hand side as $F_K(\chi)$. Then,
for $ k=0, F_K= 0$, an application of  the Lemma  ~\ref{lemma:3.1} and
$$
\inf_{\chi \in V^0_h}\|\eta - \chi\|_1=\inf_{\chi \in V^0_h}\|u - \chi\|_1\le C(q,0)h^{q-1}
$$
shows the result  for $k=0.$

For the general case, we resort to induction on $k$. Let the assertion of the theorem  be  true for $k-1$. Then, we claim that conclusion also holds   for $k$.
For all $\rho \in H^s(\Omega)\cap H^1_0(\Omega)$, we arrive from (\ref{eq:estimate-F}) at
\begin{eqnarray*}
|F_{ik}(u,\nabla u;\frac{\partial^i \eta}{\partial t^{i}},\rho)| &\le& C(max(s,0)+2,k)\sum^{k-1}_{i=0}\|\frac{\partial^i\eta}{\partial t^i}\|_{-s}\|\rho\|_{s+2}\\
&\le& C(max(q,s+2),k)h^{s+q}\|\rho\|_{s+2},\;\;\; s= -1,0,1,\ldots,r-1.
\end{eqnarray*}
Here, we have used the induction hypothesis to derive the second inequality.
Since,
$$
\inf_{\chi \in V^0_h}\|\frac{\partial^k \eta}{\partial t^k}- \chi\|_1 = \inf_{\chi \in V^0_h}\|\frac{\partial^k u}{\partial t^k}- \chi\|_1 \le C(q,k)h^{q-1},
$$
the use of lemma ~\ref{lemma:3.1} completes the  rest of the proof. \hfill{$\Box$}

Moreover, following standard argument for linear elliptic problems, see Brenner and Scott \cite{BS-2008}, the following $L^p$ for $1\leq  p \leq \infty$ estimate holds for both $\eta$ and $\eta_t$ and for $1\leq q \leq r+1$
\begin{eqnarray}\label{estimate:Lp}
\|\eta\|_{L^p}+ h \|\nabla\eta\|_{L^p} \leq C\;h^{q}\; \|u(t)\|_{W^{r+1,p}}  \;\;{\mbox {and}}\; \|\eta_t(t)\|_{L^p} \leq C\;h^{r+1}\; \|u_t(t)\|_{W^{r+1,p}}.
\end{eqnarray}

\section{Global Error Estimates}
\setcounter{equation}{0}
This section focusses on the optimal error estimates of  $u-U$ in $L^{\infty}(L^2)$ and $L^{\infty}(H^1)$-norms.

The following two Lemma  shows both $ L^\infty(L^2) $ and $L^{\infty}(H^1)$ error estimation for $ \zeta = U-\tilde{u}_h.$
\begin{lemma}\label{lemma:zeta-L2}
Let $ 2 \le q \le r+1 $. Then, 
the following estimate
\begin{eqnarray}
\|\zeta\|_{L^\infty(L^2)}+ \alpha_0 \|\zeta\|_{L^2 (H^1)} \le C(max(q,2k+1),k)\; \Big(h^q + h^{2(q-1)}\Big)
\end{eqnarray}
holds.
\end{lemma}

\noindent
{\it Proof}. From (\ref{eq:zeta}), we now rewrite it as
\begin{eqnarray}\label{eq:zeta-1}
(\zeta_t,v)+ a_{\lambda} (\tilde{u}_h,U; v)
&=&(\eta_t,v)  + \lambda (\zeta,v) -\lambda (\eta,v) \nonumber\\
&&\;\;\;+({\mathbb R}_1(\eta,\nabla \eta),\nabla v) +  ({\mathbb R}_2(\eta,\nabla \eta), v).
\end{eqnarray}
Setting $v=\zeta$ in (\ref{eq:zeta-1}), a use of coercivity  of $a_{\lambda}(\tilde{u}_h,U; \zeta)$ with the Cauchy-Schwartz inequality yields
\begin{eqnarray}\label{eq:zeta-2}
\frac{1}{2} \frac{d}{dt} \|\zeta(t)\|^2+ \alpha_0 \;\|\nabla \zeta(t)\|^2
&\leq & \Big( \|\eta_t \|_{-1} + \lambda \|\eta\|_{-1}\Big) \;\|\nabla \zeta\| + \lambda \|\zeta\|^2 \nonumber\\
&&\;\;\;+ |({\mathbb R}_1(\eta,\nabla \eta),\nabla \zeta)| + |({\mathbb R}_2(\eta,\nabla \eta), \zeta)|.
\end{eqnarray}
For the last  two terms on the right hand side of (\ref{eq:zeta-2}),  we use the form of ${\mathbb R}_1$ and ${\mathbb R}_2$ and generalized H\"older inequality with $L^4$ estimate (\ref{estimate:Lp})  to obtain
\begin{eqnarray}\label{estimate-R-1}
|({\mathbb R}_1(\eta,\nabla \eta),\nabla \zeta)| &\leq&  C \Big( \|\eta\|_{L^4}^2 + \|\nabla \eta\|_{L^4}^2\Big)\;\|\nabla\zeta\| \nonumber\\
&\leq&  C \; h^{2(q-1)} 
\;\|\nabla\zeta\|,
\end{eqnarray}
and similarly, we arrive at
\begin{eqnarray}\label{estimate-R-2}
|({\mathbb R}_2(\eta,\nabla \eta),\zeta)| 
\leq   C \; h^{2(q-1)} 
\;\|\zeta\|.
\end{eqnarray}
On substitution in (\ref{eq:zeta-2}), a use  of Young's inequality with kickback argument and  application of Gronwall's Lemma completes the rest of the proof.  \hfill{$\Box$}

Our next Lemma concerns with estimate  of $\zeta$ in gradient norm.
\begin{lemma}\label{lemma:zeta-gradient}
Let $ \max((1+d/2),2)\; \le q \le r+1$. Then, 
the following estimate holds
\begin{eqnarray}\label{eq:zeta-gradient}
&&\|\zeta_t\|_{L^2(L^2)}+ \alpha_0 \|\zeta\|_{L^\infty(H^1)}
\leq \;C(max(q,2k+2),k) \;\;h^{q}.
\end{eqnarray}
\end{lemma}
\noindent
{\it Proof}. Choosing  $v =  \zeta_t$ in (\ref{eq:zeta}),  we arrive at
\begin{eqnarray}\label{eq:zeta-t}
\|\zeta_t\|^2+ B (\tilde{u}_h,\nabla \tilde{u}_h; \zeta, \zeta_t )
&=&(\eta_t,\zeta_t)  + \lambda (\zeta, \zeta_t) -\lambda (\eta, \zeta_t) \nonumber\\
&+& ({\mathbb R}_1( e,\nabla e),\nabla \zeta_t) +  ({\mathbb R}_2(e,\nabla e), \zeta_t).
\end{eqnarray}
Note that from (\ref{eq:B-time-derivative}) with $k=1$, we obtain
\begin{eqnarray} \label{eq:B-1}
B (\tilde{u}_h,\nabla \tilde{u}_h; \zeta, \zeta_t ):= \frac{1}{2} \frac{d}{dt} B (\tilde{u}_h,\nabla \tilde{u}_h; \zeta, \zeta)- F_{01}(\tilde{u}_h, \nabla\tilde{u}_h; \zeta, \zeta).
\end{eqnarray}
On substitution (\ref{eq:B-1}) in (\ref{eq:zeta-t}), it follows that
\begin{eqnarray}\label{eq:zeta-t-1}
\|\zeta_t\|^2+ \frac{1}{2} \frac{d}{dt} B (\tilde{u}_h,\nabla \tilde{u}_h; \zeta, \zeta)
&=&\Big((\eta_t,\zeta_t)  + \lambda (\zeta, \zeta_t) -\lambda (\eta, \zeta_t)\Big) + F_{01}(\tilde{u}_h, \nabla\tilde{u}_h;\zeta, \zeta) \nonumber\\
&+& ({\mathbb R}_1(e,\nabla e),\nabla \zeta_t) +  ({\mathbb R}_2(e,\nabla e), \zeta_t) \nonumber\\
&=& I_1 (\zeta_t)+ I_2 (\zeta_t)+ I_3 (\zeta_t)+ I_4(\zeta_t).
\end{eqnarray}
For the first term on the right hand side of (\ref{eq:zeta-t-1}), apply the Cauchy-Schwarz inequality, estimates of $\eta$ and $\eta_t$ with estimate from Lemma ~\ref{lemma:zeta-L2}  to find that 
\begin{eqnarray}\label{estimate:I-1}
  |I_1(\zeta_t)| &\leq &    \Big( \|\eta_t\|+ \lambda \|\eta\| +\lambda \|\zeta\|\Big)\; \|\zeta_t\|\nonumber\\
  &\leq& C h^{2q} + \frac{1}{8}\; \|\zeta_t\|^2.
\end{eqnarray}
For the second  term on the right hand side of (\ref{eq:zeta-t-1}), a use of (\ref{eq:estimate-F}) with $i=0,k=1$ and the stability of the elliptic projection in $H^1$ yields
\begin{eqnarray}\label{estimate:I-2}
  |I_2(\zeta_t)| \leq     C(1,1)\;  \|\nabla \zeta\|^2.
\end{eqnarray}
In order to estimate the fourth  term on the right hand side of (\ref{eq:zeta-t-1}), we observe using definition of
${\mathbb R}_2 $, the generalized H\"older's inequality, the Sobolev inequality, inverse inequality  and the Young's inequality that 
\begin{eqnarray}\label{estimate:I-4}
  |I_4(\zeta_t)| &\leq&   C\; \int_{\Omega} \Big( |e|^2 + |e|\;|\nabla e| + |\nabla e|^2\Big)\; |\zeta_t|\; dx\nonumber\\
 &\leq&  C\:\Big( \|e\|_{L^4}^2 + \|\nabla e\|^2_{L^4} \Big)\; \|\zeta_t\| \nonumber\\
 &\leq& C\; (h^{2q} +h^{-d/2}\|\zeta\|^2 + h^{2(q-1)} + h^{-d/2} \|\nabla \zeta\|^2)\;  \|\zeta_t\| \nonumber\\
 &\leq& C \;\Big( h^{4(q-1)-d} + h^{-d}  \|\nabla \zeta\|^4\Big) + \frac{1}{8} \|\zeta_t\|^2.
\end{eqnarray}
For the estimate of $I_3(\zeta_t),$ rewrite it as
\begin{eqnarray}\label{estimate:I-3}
  I_3(\zeta_t)&=& \frac{d}{dt}  ({\mathbb R}_1(e,\nabla e),\nabla \zeta)  - (\frac{d}{dt}{\mathbb R}_1(e,\nabla e),\nabla \zeta).
\end{eqnarray}
Note that using the definition of $ {\mathbb R}_1(e,\nabla e)$ and the  chain rule, an application of  generalized H\"older's inequality with inverse inequality  and maximum norm bounds of $\eta$ and $\nabla \eta$ yields
\begin{eqnarray}\label{eq:I-3-1}
-(\frac{d}{dt}{\mathbb R}_1(e,\nabla e),\nabla \zeta) :=  -({\mathbb R}_1(e,\nabla e)',\nabla \zeta)
- ({\mathbb R}_{1,t}(e,\nabla e),\nabla \zeta),
\end{eqnarray}
where $ {\mathbb R}_1(e,\nabla e)'$ is the time derivative of the variable $e$ and $\nabla e$  using chain rule and ${\mathbb R}_{1,t}(e,\nabla e)$ is the time derivative of the coefficients again using chain rule.
For the estimate of the first term on the right hand side of \eqref{eq:I-3-1}, we obtain
\begin{eqnarray}\label{eq:I-3-2}
 -({\mathbb R}_1(e,\nabla e)',\nabla \zeta) &\leq& C\;\int_{\Omega} (|e_t| + |\nabla e_t|) ( |e|+ |\nabla e|)\;|\nabla \zeta|\;dx\nonumber\\
 &\leq& C\; \|\nabla e_t\|\; \|\nabla e\|_{L^{\infty}} \; \|\nabla \zeta\| \nonumber\\
  &\leq& C\;\Big( h^{q-1} + h^{-1}  \|\zeta_t\| \Big) \; \big(h^{q-1} + h^{-d/2} \|\nabla \zeta\| \big)\; \|\nabla \zeta\| \nonumber\\
  &\leq& C\;\Big( \big( h^{4(q-1)} + h^{2q}\big) + \|\nabla \zeta\|^2 + h^{-2(1+d/2)} \;\|\nabla \zeta\|^4\Big) + \frac{1}{8} \|\zeta_t\|^2.
\end{eqnarray}
For  the second term on the right hand side of \eqref{eq:I-3-1}, we  bound it as
\begin{eqnarray}\label{eq:I-3-3}
 -({\mathbb R}_{1,t}(e,\nabla e),\nabla \zeta) &\leq& C\;\int_{\Omega} (1+|e_t| + |\nabla e_t|) ( |e|^2+ |e| |\nabla e| +|\nabla e|^2)\;|\nabla \zeta|\;dx\nonumber\\
 &\leq& C\; \Big( \|e\|^2_{L^4} + \|\nabla e\|^2_{L^4}\Big)\; \|\nabla \zeta\| + C\; \Big( \|e_t\| + \|\nabla e_t\|\Big) \;\big( \|e\|^2_{L^{\infty}} + \|\nabla e\|_{L^{\infty}}\big)\;  \|\nabla \zeta\|
 \nonumber\\
 &\leq& C\;  ( \|\nabla e\|^2 + \|\nabla e\|^2_{L^4} )\; \|\nabla \zeta\|
 +  C\; \|\nabla e_t\| \; \|\nabla e\|_{L^{\infty}}\big)\;  \|\nabla \zeta\|
 \nonumber\\
&\leq& C\;\Big( h^{2(q-1)} +   (1+h^{-d/2}) \|\nabla\zeta \|^2\Big)\;   \|\nabla \zeta\| \nonumber\\
&+&  C\; \big( h^{q-1} +    h^{-1}  \|\zeta_t\| \big) \; \big(h^{q-1} + h^{-d/2} \|\nabla \zeta\| \big)\; \|\nabla \zeta\| \nonumber\\
  &\leq& C \;\Big( \big( h^{4(q-1)} + h^{2q}\big) + \|\nabla \zeta\|^2 + h^{-2(1+d/2)} \;\|\nabla \zeta\|^4\Big) + \frac{1}{8} \|\zeta_t\|^2.
\end{eqnarray}

On substitution of the estimates (\ref{estimate:I-1})-(\ref{estimate:I-3}) in  (\ref{eq:zeta-t-1}), an integration with respect to time with $\zeta(0)=0$, (\ref{estimate-R-1}), the Young's inequality, kickback arguments  and estimate of $\eta$ and its gradient yields
\begin{eqnarray}\label{eq:zeta-t-2}
 \int_{0}^{t} \|\zeta_t(s)\|^2 \;ds +  \alpha_0 \|\nabla \zeta (t)\|^2
&\leq & C(max(q,2k+1),k)\;\Big( h^{2 q} +
\int_{0}^{t}  \Big(h^{-2(1+d/2)}  \|\nabla\zeta\|^4 \nonumber\\
&+& (1+ h^{4(q-1)-2})\; \|\nabla \zeta\|^2\Big)\;ds\Big).\nonumber\\
\end{eqnarray}
An application of Gronwall's inequality with use of $L^2(L^2)$ estimate of $\nabla\zeta$ and for $(1+d/2 )\leq q \leq r+1$
shows the result and this completes the rest of the proof. \hfill{$\Box$}

A use of triangle inequality with estimates of Lemmas ~\ref{lemma:zeta-L2}-\ref{lemma:zeta-gradient} and Theorem ~\ref{theorem:eta-t} yields easily the following result.
\begin{theorem}
Let $e = U-u$, where u is the solution of (3.1) and U is the solution of (3.2) with $U(0)$ defined as in (4.3), for some k such that $\max((1+d/2),2) \le q \le r+1$
\begin{eqnarray}
\|e\|_{L^\infty(L^2)}+h \Big(\|e\|_{L^2(H^1)}+ \|e\|_{L^\infty(H^1)}\Big)
\le C(max(q,2k+2),max(k,1))\;h^q.
\end{eqnarray}
\end{theorem}
\begin{remark}
Note that  using Lemma ~\ref{lemma:zeta-L2} and Theorem ~\ref{theorem:eta-t} with triangle inequality,  we obtain
optimal in  $\|U-u\|_{L^{\infty}(L^2)}$-  estimate, when $2\leq q \leq r+1$, that is, for  $r\geq 1.$
\end{remark}
Now, we obtain the following corollary.
\begin{corollary}
 Under  the assumptions of  Lemma ~\ref{lemma:zeta-gradient}, there holds for $\max((1+d/2),2) \le q \le r+1$
\begin{equation}\label{eq:R-1-R-2}
\|{\mathbb R}_2\|_{L^2(L^2)}+ h \|{\mathbb R}_1\|_{L^2(L^2)} \le C(max(q,2k+2)) \; \Big(  h^{2(q-1)}
+ h^{(2q-d/2)} \Big).
\end{equation}
\end{corollary}

\noindent
{\it Proof}. From the definition of  ${\mathbb R}_2$,  a use of estimate of $\eta$ from
\eqref{estimate:Lp} with  \eqref{eq:zeta-gradient} and inverse estimate  yields
\begin{eqnarray}\label{eq:R-2-L-2}
\|{\mathbb R}_2\|_{L^2(L^2)} = \sup_{\|v\|_{L^2(L^2)}=1}| \int_{0}^{T} ( {\mathbb R}_2,v)\;ds|  &\leq&  C\; \Big(\int_{0}^T \int_{\Omega} \big( |e|^4 + |e|^2\;|\nabla e|^2 + |\nabla e|^4\Big)\; dx\;ds \Big)^{1/2} \nonumber\\
 &\leq&  C\:\int_{0}^{T}\big( \|e\|_{L^4}^2 + \|\nabla e\|^2_{L^4} \big)\;ds  \nonumber\\
 &\leq& C \; \int_{0}^{T}\Big ( h^{2q} + \|\zeta\|^2_{L^4} + h^{2(q-1)} + h^{-d/2} \|\nabla \zeta\|^2\Big)\;ds
 \nonumber\\
 &\leq& C(T)  \; \Big ( h^{2q} + \|\nabla\zeta\|^2_{L^2(L^2)} + h^{2(q-1)} + h^{-d/2} \|\nabla \zeta\|^2_{L^2(L^2)} \Big) \nonumber\\
&\leq& C \:\Big ( h^{2(q-1) } + h^{2q-d/2} \Big).
\end{eqnarray}
Similarly,  the estimate
\begin{equation} \label{eq:R-1-H-1}
\|{\mathbb R}_1\|_{L^2(H^{-1})} \leq C \:\Big ( h^{2(q-1) } + h^{2q-d/2} \Big)
\end{equation}
holds and the rest of the proof follows.
 \hfill{$\Box$}

As a consequence of superconvergence result of $\|\nabla\zeta\|$ from \eqref{eq:zeta-gradient},  the following maximum norm estimate is derived which we put as a Theorem.
\begin{theorem} There holds for  $2\leq q\leq r+1$
\begin{equation}\label{eq:e-max-norm}
\|U-u\|_{L^{\infty}(L^{\infty})} 
\leq C(max(q,2k+2),max(k,1)) \Big(\log(1/h)\Big)^{m}\;\;h^q,
\end{equation}
where $m=0$ for $d=1$,  and $m=1$ when $d=2.$
\end{theorem}

\noindent
{\bf Proof}. From the superconvergence result \eqref{eq:zeta-gradient},
we obtain for $d=1$, the maximum norm estimate of $\zeta$, that is,
\begin{equation}\label{eq:zeta-max-norm}
\|\zeta\|_{L^{\infty}(L^{\infty})} \leq \|\zeta\|_{L^{\infty}(H^1)} \leq C(max(q,2k+2),max(k,1))\;h^q.
\end{equation}
Therefore, a use of triangle inequality with \eqref{eq:zeta-max-norm} and \eqref{estimate:Lp} for $p=\infty$ yields for 1D-problem the following max norm estimate of $e$ for $2\leq q\leq r+1$
\begin{equation}\label{eq:e-max-norm}
\|U-u\|_{L^{\infty}(L^{\infty})} 
\leq C(max(q,2k+2),max(k,1))\;h^q.
\end{equation}

For $d=2$, using Sobolev inequality, we arrive from the superconvergence result \eqref{eq:zeta-gradient}, the max norm estimate
\begin{equation}\label{eq:zeta-max-norm-2D}
\|\zeta\|_{L^{\infty}(L^{\infty})} \leq C \log(1/h)\;\|\zeta\|_{L^{\infty}(H^1)} \leq C(max(q,2k+2),max(k,1))\;
 \log(1/h)\;\;h^q.
\end{equation}
Hence, we derive the max norm estimate for 2D-problem for $2\leq q\leq r+1$ as
\begin{equation}\label{eq:e-max-norm-2D}
\|U-u\|_{L^{\infty}(L^{\infty})} 
\leq C(max(q,2k+2),max(k,1)) \log(1/h)\;\;h^q.
\end{equation}
This completes the rest of the proof. \hfill{$\Box$}

\section {The Quasi-projection}
\setcounter{equation}{0}
Let  $z_0= \eta$ and $\theta_0 = \zeta$. Define maps $z_j : J \rightarrow V^0_h$ recursively by
\begin{eqnarray}\label{eq:z-j}
B(u,\nabla u; z_j,v) = -(\frac{\partial z_{j-1}}{\partial t},v), \;\; v \in V^0_h,\;\; t\in J,\;\; j = 1,2,\ldots.
\end{eqnarray}
\begin{theorem}\label{theorem:z-j}
Let $j \ge 0, k\ge 0,1 \le q \le r+1$ and assume that $\frac{\partial^{j+k} u}{\partial t^{j+k}} \in H^q(\Omega)$ for $t \in J$. Then for $-1 \le s \le r-1-2j$, the following estimate holds
\begin{eqnarray}\label{eq:estimate-z-j}
\|\frac{\partial^k z_j}{\partial t^k}\|_{-s}= C\Big(max(q,max(s,0)+2j+2),k+j \Big)\;h^{s+q+2j}.
\end{eqnarray}
\end{theorem}

\noindent
{\it Proof}. The proof is carried out by induction $j$. For $j=0,z_0 = \eta$ and this case is covered by Theorem 3.1. Now suppose that (\ref{eq:estimate-z-j}) is true for $j-1$, then we show (\ref{eq:estimate-z-j}) to be true for $j$.
For $j >0$ and $ k \ge 0$, let
$$
F(\rho)= - \sum_{i=0}^{k-1} F_{ik}(\frac{\partial^i z_j}{\partial t^i},\rho)-(\frac {d^k}{dt^k}(\frac{\partial z_{j-1}}{\partial t}),\rho), \quad \rho \in H^{s+2}(\Omega)\cap H^1_0(\Omega).
$$
Then,
$$
B(u,\nabla u;\frac{\partial^k z_j}{\partial t^k},\chi) = F(\chi),\quad \chi \in V^0_h + H^{s+2}(\Omega).
$$
Now
\begin{eqnarray*}
|F(\rho)|&\le& | \sum_{i=0}^{k-1} F_{ik}(\frac{\partial^i z_j}{\partial t^i},\rho)|+|(\frac {d^k}{dt^k}(\frac{\partial z_{j-1}}{\partial t}),\rho)|\\
&\le& C(max(s,0)+2,k)\sum_{i=0}^{k-1}\|\frac{\partial^i z_j}{\partial t^i}\|_{-s}\|\rho\|_{s+2}+|(\frac{\partial^{k+1}}{\partial t^{k+1}} z_{j-1},\rho)|\\
&\le& \{C(max(s,0)+2,k)\sum_{i=0}^{k-1}\|\frac{\partial^i z_j}{\partial t^i}\|_{-s}+\|\frac{\partial^{k+1}}{\partial t^{k+1}} z_{j-1}\|_{-s-2}\}\|\rho\|_{s+2}\\
&\le& \{C(max(s,0)+2,k)\sum_{i=0}^{k-1}\|\frac{\partial^i z_j}{\partial t^j}\|_{-s} +C(max(q,max(s,0)+2j+2),k+j)h^{s+q+2j} \} \|\rho\|_{s+2}.
\end{eqnarray*}
The  last inequality uses induction hypothesis.
If we consider $ p= r-2j$
$$
M_{s+2}=C(max(q,max(s,0)+2j+2),k+j)(h^{s+q+2j} + \sum_{i=0}^{k-1}\|\frac{\partial^i z_j}{\partial t^i}\|_{-s})
$$
for $s=-1,0,1,\ldots,r-2j-1,$
in the Lemma 3.1, then $F$ fulfills all the hypothesis. Therefore, we obtain
$$
\|\frac{\partial^k z_j}{\partial t^k}\|_{-s} \le C(max(s,0)+2,0)\left[(M_1+\inf_{\chi \in V^0_h}\|\frac{\partial^k z_j}{\partial t^k}-\chi\|)h^{s+1}+M_{s+2}\right]
$$
for $ s= -1,0,\ldots,r-2j-1$. As $z_j \in V^0_1 $, the infimum appearing here is zero and it follows that
\begin{eqnarray*}
\|\frac{\partial^k z_j}{\partial t^k}\|_{-s} &\le& C(max(q,max(s,0)+2j+2),k+j)[h^{s+q+2j} \\
&+&\sum_{i=0}^{k-1}\|\frac{\partial^i z_j}{\partial t^i}\|_1 h^{s+1}+\sum_{i=0}^{k-1}\|\frac{\partial^i z_j}{\partial t^i}\|_{-s}].
\end{eqnarray*}
For $k = 0$, the inequality (\ref{eq:estimate-z-j}) follows. For $k \ge 1$, the theorem can be completed by simple induction on $k$. This completes the induction on $j$ and hence, the inequality (\ref{eq:estimate-z-j}) is proved.
This completes the rest of the proof.  \hfill{$\Box$}

For one of the main superconvergent result, we can choose the initial condition $U$ at $t=0$ as
\begin{eqnarray}\label{eq:initial-condition}
U(0)=\tilde u(0) + z_1(0) + \ldots + z_k(0), \;\;\;{\mbox {for}}\; 2k \le r-1.
\end{eqnarray}
Let $\theta_k = \zeta +z_1+\ldots +z_k$. Then for all $v \in V_h^0$,
\begin{eqnarray}\label{eq:theta-k}
(\frac{\partial \theta_k}{\partial t},v)+B(u,\nabla u;\theta_k,v)&=&(\frac{\partial z_k}{\partial t},v)-\lambda (z_k,v) + \lambda (\theta_k,v)\nonumber\\
&&\;\;\;+({\mathbb R}_1(e,\nabla e),\nabla v)+({\mathbb R}_2(e,\nabla e, v),
\end{eqnarray}
where $$({\mathbb R}_1(e,\nabla e), \nabla v)+({\mathbb R}_2(e,\nabla e),v) = -\displaystyle{\int_\Omega[\nabla v^T,v]\int_0^1(1-s){\mathbb A}^{'}(w,\nabla w)ds
\left [\begin{array}{c}
\nabla e\\
e \end{array}\right ]
\left [\begin{array}{c}
\nabla e\\
e \end{array}\right ] \;dx} $$
with $\theta_k(0)=0.$
The above relation follows as directed consequence of (\ref{eq:zeta}) and the definition of the $z_j$'s.

The following theorem gives a bound for $\theta_k$.
\begin{theorem}\label{theorem:estimate-theta-k}
Let $ 2k \le r-1$ and $ 1 \le q \le r+1$. 
Then, the following estimate
\begin{eqnarray} \label{eq:estimate-theta-k}
\|\theta_k\|_{L^{\infty}(L^2)}+\alpha_0\; \|\theta_k\|_{L^2 (H^1)} &\le&  C(max(q,2k+3),k+1)h^{q+min(2k+1,r-1)}\nonumber\\
&+& C \; \Big ( h^{2(q-1) } + h^{2q-d/2} \Big)
\end{eqnarray}
holds.
\end{theorem}

\noindent
{\it Proof}. Choose $v = \theta_k$ in (\ref{eq:theta-k}) and then, we obtain
\begin{eqnarray*}
\frac{1}{2}\frac{d}{dt}\|\theta_k\|^2 +B(u,u;\theta_k,\theta_k)&=& (\frac{\partial}{\partial t}z_k,\theta_k) + \lambda \|\theta_k\|^2 - \lambda (z_k, \theta_k) \\
&&\;\;+({\mathbb R}_1(e,\nabla e),\nabla \theta_k)+({\mathbb R}_2(e,\nabla e),\theta_k).
\end{eqnarray*}
A use of the coercivity of  $B (u,\nabla u;\cdot,\cdot)$ with the Youngs inequality 
and the kickback argument yields
$$
\frac{d}{dt}\|\theta_k\|^2 +\alpha_0\; \|\theta_k\|_1^2 \le C \;\Big(\|\frac{\partial}{\partial t}z_k\|^2_{-1}+ \|{\mathbb R}_1(e,\nabla e)\|_{L^2}^2 + \|{\mathbb R}_2(e,\nabla e)\|_{H^{-1}}^2\Big) + C\; \|\theta_k\|^2.
$$
Integrating with respect to $t$ and using the estimate (\ref{eq:estimate-z-j}) for $k=1,j=k,s=1,$ an application of Gronwall's inequality with estimates \eqref{eq:R-2-L-2} and \eqref{eq:R-1-H-1} yields the desired results. This concludes the rest of the proof.
\hfill{$\Box$}
\begin{remark}
From the above theorem that is, Theorem ~\ref{theorem:estimate-theta-k}, it is observed  for $d=1,2, 3$ with $\min(2k+1, r-1)$ as $r-1$  and $q=r+1$ that super-convergence result $\|\theta_k(t)\| =O (h^{2r})$ holds for $r\geq 2.$
\end{remark}
As a consequence, error estimates in Sobolev spaces of negative index are easily derive, which are given in the form of a corollary.  
\begin{corollary} \label{cor:negative-norm}
Let $2k\leq r-1$, $1\leq q \leq  r+1$ and $U(0)$ be defined by \eqref{eq:initial-condition}. Then, there holds for $0\leq s\leq \min(2k+1,r-1)$
\begin{equation}\label{eq:negative-norm}
\|(U-u)\|_{L^{\infty}(H^{-s}(\Omega))} \leq \;C(max(q,2k+3),k+1)h^{q+s}.
\end{equation}
\end{corollary}

\noindent
{\it Proof}. Since
$$U-u :=-\theta_k +  \eta +\sum_{j=1}^{k} z_j,$$
a use of estimate \eqref{eq:estimate-theta-k} with estimates  \eqref{eq:eta-negative} and  \eqref{eq:estimate-z-j} concludes the result. This completes the proof. \hfill{$\Box$}

\section{Super Convergence Result for the case of Single Space Variable.}
\setcounter{equation}{0}
Consider the case $\Omega = (0,1):=I$ with the finite element subspace being piecewise polynomial functions of degree $r$. Let $\Pi_h = \{x_0,x_1,\ldots,x_{N_h}\}, 0=x_0 <x_1 < \ldots < x_{N_h}=1$ with $h_i = x_i -x_{i-1}$ and $\max_{1 \le i \le N_h} h_i=h$. Further, assume that the partition is quasi-uniform, that is, $\max_{1 \le i,j \le N_h}( \frac{h_j}{h_i}) \le C,$ where the constant $C$ is independent of $h$. Assume that $V^0_h$ consists of $H_0^1(\Omega)$ functions of $I$ whose restrictions on $I_j$'s are polynomials of degree atmost $r$, where $I_j=(x_{j-1},x_j)$. At each knot $x_i,1 \le i \le N_h -1 ,$ the element of $V_h^0$ will be assume to be $C^{p_i}$ functions, where $0 \le p_{i} \le r $, and $ p_0 =p_{N_h}=0$. However, at any knot at which the superconvergence is to take place, the smoothness constraint on $V_h^0$ must be restricted to continuity only. The properties \eqref{eq:approx-property} and \eqref{eq:inverse}  hold  true for a such a choice of $V_h^0$.

Following  the procedure, developed by Douglas et al. \cite{DDW-1978},  we now establish the  knot superconvergence  result  of this section.

Let $\bar{x}\in I$ be a knot in each of the partition, that is, for each $h$, there is $i(h)$ such that
 $\bar x = x_{i(h)}.$  Further, let $p_{i(h)}=0$. Denote the space $\tilde H^s$ by
\begin{eqnarray}\label{eq:tilde-H-s}
\{u:u|_{(0,\bar x)} \in H^s((0,\bar x)),u|_{(\bar x,1)} \in H^s((\bar x,1))\}\times{\mathbb R}
\end{eqnarray}
and it's norm be given by
\begin{eqnarray} \label{eq:triple-norm}
|\|(u,\beta)\||^2 = \|u\|^2_{H^s((0,\bar x))}+\|u\|^2_{H^s((\bar x,1))}+\beta^2.
\end{eqnarray}
For any element $(u,\beta)$ and $(v,\gamma)\in \tilde H^0$, define the inner product by
\begin{eqnarray}\label{eq:norm-0}
\left [ (u,\beta),(v,\gamma)\right] = (u,v)+\beta\gamma,
\end{eqnarray}
where $(u,v)$ denotes $L^2$ inner product.

The space $\tilde H^{-s},$  for $s \ge 0$ is defined by duality with respect to the above inner product, see \cite{AD-1979}. For $s \ge 0$ and $z \in H^1(I)$, the norm on $\tilde H^{-s}$ is given by
\begin{eqnarray}\label{eq:z-triple-norm-negative}
|\|z\||_{-s} = \sup \{ \frac{[(z,z(\bar x)),(u,\beta)]}{|\|(u,\beta)\||_s}:|\|(u,\beta)\||_s\not= 0\}
\end{eqnarray}
It is easy to see
\begin{eqnarray}\label{eq:z-bar}
|z(\bar x)| \le |\|z\||_{-s} , \quad s \ge 0.
\end{eqnarray}

Now we derive an estimate analogous to the one in Theorem  ~\ref{theorem:z-j} for the triple norm.
\begin{theorem}\label{theorem:z-j-s}
Let $ 1 \le q \le r+1$ and $0\le s \le r-2j-1$. Then for $ 0<h<\epsilon$ sufficiently small and $t \in J $,
\begin{eqnarray}\label{eq:z-j-s}
|\|\frac{\partial^k z_j}{\partial t^k}\||_{-s} \le C(max(q,s+2j+2),j+k)h^{s+q+2j}
\end{eqnarray}
where $z_j$ is defined by \eqref{eq:z-j}.
\end{theorem}

\noindent
{\it  Proof}. For $(\psi,\beta)\in \tilde H^s$, define $\phi \in H_0^1$ by
\begin{eqnarray}\label{eq:adjoint}
L^{*}(u)\phi=\psi,x \in I \setminus {\bar x},\;\;\;\;\;\frac {\partial A}{\partial \xi}(u,u_x)\phi_x|^{\bar x+0}_{\bar x-0}=-\beta.
\end{eqnarray}
Then, $B(u,u_x;\mu,\phi)=[(\mu,\mu(\bar x)),(\psi,\beta)] $, for $ \mu \in H^1_0(I)$.\\
By regularity theorem, it follows that
\begin{eqnarray} \label{eq:regularity}
|\|(\phi,\phi(\bar x))\||_{s + 2} \le C(s+2,0)|\|(\psi,\beta)\||_s
\end{eqnarray}
Note that $\rho \in W^{s+1,\infty}(I),\;\; \mu \in H_0^1(I) $
\begin{eqnarray*}
(\rho\mu_x,\phi_x) 
&=& \rho(\bar x)\mu(\bar x)\phi_x|^{\bar x+0}_{\bar x-0} -\int_0^{\bar x}\mu(\rho\phi_x)_x \;dx
- \int_{\bar x}^1\mu(\rho\phi_x)_x\;dx\\
&=& \frac{1}{A_\xi (u(\bar x),u_x(\bar x))}\rho(\bar x)\mu(\bar x)\beta -\int_0^{\bar x}\mu(\rho\phi_x)\; dx
-\int_{\bar x}^1\mu(\rho\phi_x)_x\;dx
\end{eqnarray*}
and hence,
\begin{eqnarray}
\label{eq:6.9}
|(\rho\mu_x,\phi_x)| &\le& C|\rho(\bar x)\mu(\bar x)\beta|+|(\mu,(\rho\phi_x)_x)|
 \nonumber \\
&\le& |\|\mu\||_{-s}|\|(\psi,\beta)\||_s
\end{eqnarray}
where the constant C depends upon $\|\rho\|_{W^{s+1,\infty}},$ but no higher derivatives of $\rho$.

 In order to  prove the theorem, we use   induction in two ways, that is, an outer induction on $ j $ and an inner induction on $k.$

For $j=0$ and $v \in V^0_h$, it follows  from \eqref{eq:adjoint}, the definition of $B$ and its time derivative  that
\begin{eqnarray*}
[(\frac{\partial^k \eta}{\partial t^k},\frac{\partial^k \eta}{\partial t^k}(\bar x)),(\psi,\beta)] &=& [(\frac{\partial^k \eta}{\partial t^k},\frac{\partial^k \eta}{\partial t^k}(\bar x)),(L^{*}(u)\phi,\beta)]\\
&=&(\frac{\partial^k \eta}{\partial t^k},L^{*}(u)\phi)+\frac{\partial^k \eta}{\partial t^k}(\bar x)\beta \\
&=& B(u,u_x,\frac{\partial^k \eta}{\partial t^k},\phi)\\
&=& B(u,u_x,\frac{\partial^k \eta}{\partial t^k},\phi - v)+ \sum_{i=1}^{k-1} F_{ik}(\frac{\partial^i \eta}{\partial t^i},\phi -v)-\sum_{i=1}^{k-1} F_{ik}(\frac{\partial^i \eta}{\partial t^i},\phi),\quad v \in V_h^0.
\end{eqnarray*}
Now using \eqref{eq:6.9} for this last term on the right hand side, we  arrive at
\begin{eqnarray}\nonumber
[(\frac{\partial^k \eta}{\partial t^k},\frac{\partial^k \eta}{\partial t^k}(\bar x)),(\psi,\beta)]&\le& C\;\|\frac{\partial^k \eta}{\partial t^k}\|_1 \;\inf_{v \in V_h^0}\;\|\phi - v\|_1 \\ \nonumber
&+&\sum_{i=1}^{k-1} C(1,k-1)\;\|\frac{\partial^i \eta}{\partial t^i}\|_1\inf_{v \in V^0_h }\|\phi-v\|_1 \\
& +&\sum_{i=1}^{k-1} C(s+2,k-i)\;|\|\frac{\partial^i \eta}{\partial t^i}\||_{-s}\;|\|(\psi,\beta)\||_s.\label{eq:estimate-j-0}
\end{eqnarray}
The assumption that the elements of $V^0_h$ are only continuous functions and not differentiable at $ x = \bar x$ implies that
\begin{equation*}
\inf_{v \in V_h^0} \; \|\phi- v\|_1 \leq   C\; h^{s+1} \;|\| (\phi,\phi(\bar{x})\||_{s+2},
\end{equation*}
and hence,
\begin{eqnarray} \label{eq:estimate-eta-k}
|\| \frac{\partial^k \eta}{\partial t^k}\||_{-s} 
&\le& C(s+2,0)\|\frac{\partial^k \eta}{\partial t^k}\|_1h^{s+1}+\sum_{i=1}^{k-1} C(s+2,k-i)\|\frac{\partial^i \eta}{\partial t^i}\|_1 h^{s+1}\nonumber\\
&+& \sum_{i=1}^{k-1} C(s+2,k-i)\|\frac{\partial^i \eta}{\partial t^i}\|_{-s}.
\end{eqnarray}
The case for $j=0$ follows from the Theorem ~\ref{theorem:eta-t} and induction hypothesis on $k.$

For $j > 0$, the proof is quite similar to the above that is for $v \in V_h^0$
\begin{eqnarray*}
[(\frac{\partial^k z_j}{\partial t^k},\frac{\partial^k z_j}{\partial t^k}(\bar x)),(\psi,\beta)]&=&B(u,u_x;\frac{\partial^k z_j}{\partial t^k},\phi)\\
&=&B(u,u_x;\frac{\partial^k z_j}{\partial t^k},\phi-v)+\sum_{i=1}^{k-1}F_{ik}(\frac{\partial^i z_j}{\partial t^i},\phi - v)\\
&-&\sum_{i=1}^{k-1}F_{ik}(\frac{\partial^i z_j}{\partial t^i},\phi)+(\frac{d^k}{dt^k}(\frac{\partial z_{j-1}}{\partial t}),\phi - v)\\
&-&(\frac{d^k}{dt^k}(\frac{\partial z_{j-1}}{\partial t}),\phi).
\end{eqnarray*}
Similar to the case $j=0$, now an appropriate application of \eqref{eq:adjoint} with \eqref{eq:regularity}, \eqref{eq:6.9} and \eqref{eq:estimate-j-0} yields
\begin{eqnarray}\nonumber
|\|\frac{\partial^k z_j}{\partial t^k}\||_{-s} &\le& C(s+2,0)\|\frac{\partial^k z_j}{\partial t^k}\|_1 h^{s+1} + \sum_{i=1}^{k-1}C(s+1,k-i)\|\frac{\partial^i z_j}{\partial t^i}\|_1 h^{s+1}\\ \nonumber
&+&\sum_{i=1}^{k-1}C(s+2,k-i)|\|\frac{\partial^i z_j}{\partial t^i}\||_{-s}+C(s+1,0)|\|\frac{\partial^{k+1} }{\partial t^{k+1}}z_{j-1}\||_{-1}h^{s+1}\\
&+&C(s+1,0)|\|\frac{\partial^{k+1}}{\partial t^{k+1}}z_{j-1}\||_{-s-2}.
\end{eqnarray}
Now the conclusion of the theorem follows from Theorem ~\ref{theorem:z-j} and by Induction.
\hfill{$\Box$}

From \eqref{eq:z-bar}  and \eqref{eq:z-j-s}, it follows that
\begin{eqnarray}\label{eq:estimate-z-j-bar}
|z_j(\bar x,t)| \le C(max(q,s+2j+1),j)h^{q+s+2j},\quad j=0,1,2,\ldots.
\end{eqnarray}

The next Theorem constitutes the main result of the present work.
\begin{theorem}\label{theorem:knot-super-cgt}
Let $1 \le q \le r+1$ and k be an integer satisfying $ 0 \le k \le [\frac{r-1}{2}]$. Let $u$ be a sufficiently regular solution of \eqref{eq:1.1}  and $U$ the Galerkin approximate solution is given by \eqref{eq:semi-discrete}  with $U(0)=\tilde u(0)+z_1(0)+z_2(0)+\ldots +z_k(0),$ with $z_k^{'}$ is defined as in \eqref{eq:initial-condition}. Then, there is a positive constant $C$ such that 
\begin{eqnarray}\label{eq:knot-super-cgt}
|(u-U)(\bar x,t)| \le \left\{\begin{array} {ll}
 \displaystyle C(max(q,2k+3,k+1))h^{-\frac{1}{2}}(h^{q+2k}+h^{2q-2}); &2k \le r-2   \\
\displaystyle C(max(q,2k+3,k+2))h^{-\frac{1}{2}}(h^{q+2k+1}+h^{2q-2});&2k \le r-1
\end{array} \right.
\end{eqnarray}
where $\bar x = x_i(h)$ is a knot at which the smoothness of $V^0_h$ reduces to continuity.\\
\end{theorem}

\noindent
{\it Proof}. Writing $u-U$ in the form
\begin{eqnarray}
(u-U)(\bar x, t)= \theta_k(\bar x,t)-\sum_{j=0}^k z_j(\bar x,t)
\end{eqnarray}
and using inverse property \eqref{eq:inverse}, for the finite dimensional space $V^0_h,$  we obtain
\begin{eqnarray}\label{eq:theta-bar}
|\theta_k(\bar x,t)|\le Ch^{-\frac{1}{2}}\|\theta_k(\cdot,t)\|.
\end{eqnarray}
On combining \eqref{eq:estimate-z-j-bar}, \eqref{eq:theta-bar} and \eqref{eq:estimate-theta-k}, we arrive at the required result and this concludes the rest of  the proof. \hfill{$\Box$}

\begin{corollary}
For $k =[\frac{r-1}{2}]$ and $q=r+1$
\begin{eqnarray}
|(u-U)(\bar x,t)| \le \left\{\begin{array} {ll}
 \displaystyle C(r+1,\frac{1}{2}(r+1))h^{2r-\frac{1}{2}},  & \mbox{for odd  }  r \\
\displaystyle  C(r+1,\frac{1}{2}(r+1))h^{2r-\frac{1}{2}},  & \mbox{for even  } r.
\end{array} \right.
\end{eqnarray}
\end{corollary}

\begin{remark}
There is a deterioration of order $\frac{1}{2}$ in the superconvergence as compared to the order in Arnold and Douglas \cite{AD-1979}.  However,  in stead of using  \eqref{eq:theta-bar}, one can exploit the superconvergence property of $\|\theta_k\|_1$ norm, but we still loose half power in the superconvergence result. Therefore, we believe that this is due to the gradient term present in the operator $A$ in the equation \eqref{eq:1.1}.
\end{remark}

\section{Conclusion}
Based on Taylor's expansion, an elliptic projection is developed for the steady state linearized problem and related error estimates are derived. Optimal error estimates in both $L^{\infty}(L^2)$ and $L^{\infty}(H^1)$-norms are proved  for the quasilinear parabolic problems with nonlinearity depending on gradient by employing   the Gronwall  type result, which hold true for $d=1,2,3$. These results improve upon the earlier results of  \cite{A-1977}. A  use of quasi-project technique yields a superconvergence result for the error between the Galerkin approximation and quasi-projection and as a consequence, negative norm estimates between the exact solution and semidiscrete Galerkin approximation are established. It is observed that the Bramble-Schatz post processing, see \cite{BS-1977}  combined with negative norm estimate may provide superconvergence results, which can be achieved provided some estimates of the difference  quotients are shown to be true for this nonlinear problem.
In a single space variable, superconvergence results at the nodal points are derived. Developing higher order time discretization method  for the completely discrete schemes combined with quasi-projection tool shall be more challenging and will be a part of our future endevour.

\vspace{2em}
\noindent
{\bf Acknowledgement}. This research of the second author is supported by Chiangmai University and the Centre of Excellence in Mathematics, The Commmision on Higher Education, Thailand.

\bibliographystyle{plain}

\end{document}